\documentclass[11pt]{article}
\usepackage{amsbsy,amsmath,amssymb,dsfont}
\begin{document}
\title{\Large \bf Equivalence of Strong Brunn-Minkowski Inequalities and  CD Conditions in Heisenberg Groups
\thanks{Supported  partially by NNSF of China }}

\author{\small \bf Juan Zhang;
\small \bf Peibiao Zhao
}
\date{}
\maketitle

\vskip 20pt

\begin{center}
\begin{minipage}{12cm}
\small
 {\bf Abstract:} The present paper investigates  the sub-Riemannian version of the equivalence between the curvature-dimension conditions 
 and    strong Brunn-Minkowski inequalities in the sub-Riemannian Heisenberg group $\mathbb{H}^n$.

 We adopt the optimal transport and approximation of $\mathbb{H}^n$ developed by Ambrosio and Rigot \cite {AR} and combine the celebrated works by M. Magnabosco,  L. Portinale and T. Rossi \cite {MPR} to confirm this.



 \vskip 0.1cm
\end{minipage}
\end{center} 

\section{\bf Introduction}

~~~~The property for a space to satisfy a so-called curvature-dimension condition CD(K, N) is interpreted as behaving in some aspects as 
a Riemannian manifold with dimension $\leq N$ and Ricci curvature $\geq K$  at any point, i.e., this Riemannian manifold have Ricci curvature bounded below and dimension bounded above. From the creativity papers by Lott and Villani \cite {LV} and Sturm \cite {S1,S2}, they defined a new notion of curvature-dimension CD(K, N) using optimal transport. The two most critical advantages of their use of this theory are: (i) the theory can be developed on very general sets (typically 
on Polish metric spaces $(X,\rho)$), (ii) the geodesics of the Wasserstein space (a metric space made of the probability measures used in optimal transport) are represented as a probability measure in the space of the geodesics of $(X,\rho)$. Up to now there are four common types of entropy: Rényi-type entropy, Shannon entropy, kinetic-type entropy and Tsallis entropy. The same type of curvature-dimension condition can be established by choosing different entropies. In the Riemannian setting \cite{CM}, metric measure spaces with generalized lower Ricci curvature bounds support various geometric and functional inequalities including Borell-Brascamp-Lieb, Brunn-Minkowski, Bishop-Gromov inequalities. It is a fundamental question whether the method used in \cite {LV,MPR,S1,S2}, based on optimal mass transportation works in the setting of singular spaces with no lower curvature bounds. A large class of such spaces are the sub-Riemannian geometric structures or Carnot-Carath\'{e}odory geometries, see Gromov \cite {G}.

Typical example of sub-Riemannian spaces are  the Heisenberg group,  Carnot group and Grushin space, etc. This paper focuses on the Heisenberg group, so let's briefly introduce the Heisenberg group. The Heisenberg group is a noncommutative stratified nilpotent Lie group. As a set it can be identified with its Lie algebra $\mathbb{R}^{2n+1}\cong\mathbb{C}^n\times\mathbb{R}$, via exponential coordinates. We denote a point in $\mathbb{H}^n$ by $x=(\xi,\eta,t)=(\zeta,t)$ where $\xi=(\xi_1,\ldots,\xi_n)$, $\eta=(\eta_1,\ldots,\eta_n)\in\mathbb{R}^n$, $t\in\mathbb{R}$ and $\zeta=(\zeta_1,\ldots,\zeta_n)\in\mathbb{C}^n$ with $\zeta_j=\xi_j+i\eta_j$. The group law in this system of coordinates is given by
$$ x\ast y =(\zeta,t)\ast(\zeta',t'):=(\zeta+\zeta',t+t'+2\sum \limits^n \limits_{j=1} Im\zeta_j \bar{\zeta'}_j ),~~~\forall~x,~y\in\mathbb{H}^n.$$
The center of the group is $C=\{(\zeta,t)\in\mathbb{H}^n: \zeta=0\}$, the neutral element of $\mathbb{H}^n$ is $0_{\mathbb{H}^n}=(0_{\mathbb{C}^n},0)$ and the inverse element of $(\zeta,t)$ is $(-\zeta,-t)$. A left translation by $x\in\mathbb{H}^n$ is the mapping $L_x: \mathbb{H}^n\rightarrow\mathbb{H}^n$, $y\rightarrow L_x(y)=x\ast y$, similarly $R_x: \mathbb{H}^n\rightarrow\mathbb{H}^n$, $y\rightarrow R_x(y)=y\ast x$ is right translation. For any $\lambda >0$, the mapping $\delta_\lambda : \mathbb{H}^n\rightarrow\mathbb{H}^n$, $x\rightarrow \delta_\lambda (\zeta,t)=(\lambda\zeta,\lambda^2 t)$ is called a dilation. For more details bibliographic references, we refer the reader to \cite {CDPT,J1,J2,M}.

During the past decade considerable effort has been made to establish some geometric and functional inequalities on sub-Riemannian spaces. The quest for Borell-Brascamp-Lieb and Brunn-Minkowski type inequalities became a hard nut to crack even on simplest sub-Riemannian setting such as the Heisenberg group $\mathbb{H}^n$ endowed with the usual Carnot-Carath\'{e}odory metric $d$ and the Haar measure $\mathcal L^{2n+1}$. (As the Heisenberg group is a non-abelian group, in this paper we have chosen the Haar measure $\mathcal L^{2n+1}$.) One of the reasons for this is that although there is a good first order Riemannian approximation (in the pointed Gromov-Hausdorff sense) of the sub-Riemannian metric structure of the Heisenberg group $\mathbb{H}^n$, there is no uniform lower bound on the Ricci curvature in these approximations (see \cite {CDPT}, Section 2.4.2); indeed, at every point of~~$\mathbb{H}^n$ there is a Ricci curvature whose limit is $-\infty$ in the Riemannian approximation. In the absence of a uniform lower Ricci curvature bound, the Riemannian Borell-Brascamp-Lieb inequality and the Brunn-Minkowski inequality are not easy to generalize directly in the context of the Heisenberg group. There are many celebrated works about such problems in the Heisenberg group, we can see \cite {J1,J2,J3,J4,M}. Recently, Balogh et al. \cite {BKS} yield good results via optimal mass transportation, by using convergence results for optimal transport maps in the Riemannian approximation of~~$\mathbb{H}^n$ due to Ambrosio and Rigot \cite {AR}, they introduce a correct sub-Riemannian geometric quantities which can replace the lower curvature bounds and can be successfully used to establish the geodesic version of Borell-Brascamp-Lieb, Pr$\acute{e}$kopa-Leindler, Brunn-Minkowski and entropy inequalities in the Heisenberg group $\mathbb{H}^n$.

In the Heisenberg group, the Brunn-Minkowski inequality is a concavity property of the mass of the s-intermediate points, namely of the function $\mathcal L^{2n+1}(Z_s(A,B))^\frac{1}{2n+1}$ ($A$, $B$ are two nonempty Borel sets of~~$\mathbb{H}^n$ and the set $Z_s(A,B)$ see (3.3)), see Theorem 3.7. A distinctive feature of the Brunn-Minkowski inequality is that its formulation does not involve optimal transport, and its proof relies on the set of $s$-intermediate points $Z_s(A, B)$ of the sets $A$ and $B$. It is well known that it contains $spt((T_s)_\sharp\eta)$, the tight support of the Wasserstein $s$-intermediate points. Therefore, a natural strengthening of the Brunn-Minkowski inequality is to require that the aforementioned concavity property holds, not for the whole set of s-midpoints, but only for the support of the Wasserstein s-intermediate point. A new result of this paper is a strong Brunn-Minkowski inequality on the tight support of Wasserstein s-intermediate points in the Heisenberg group $\mathbb{H}^n$, which we denote by SBMI, cf. Theorem 3.10. Although this SBMI still depend on optimal transport, seeing it reminds you of the Brunn-Minkowski inequality.

In this paper, our main result (Theorem 4.1) is that on Heisenberg group, the SBMI is equivalent to the CD condition in the setting of essentially non-branching metric measure spaces. The main differences between our results and the Riemannian case (see \cite {MPR}) are as follows:\\
\noindent{\bf Remark 1.1}.~(i) The Heisenberg group is a non-abelian group, and the measure we choose in this paper is Haar measure~~$\mathcal L^{2n+1}$. As the left translation map is affine and its vectorial part has the determinant 1, it follows that the Haar measure of~~$\mathbb{H}^n$ is the Lebesgue measure of~~$\mathbb{R}^{2n+1}$ which is left (and actually also right) invariant.

(ii) We need to emphasize the difference between the Riemannian manifold and Heisenberg group versions of distortion coefficients, which encodes information on the geometric structure of the space. In the Riemannian manifold, the coefficient $\tau^{K,N}_s(d(x,y))$ is a function about the distance between two points $x$ and $y$. But in the Heisenberg group, Ricci curvatures tend to $-\infty$ and $\lim_{K\rightarrow-\infty}\tau^{K,N}_s(d(x,y))=0$ for every $s\in[0,1]$, to avoid this situation, we use the angle $\theta(x,y)$ to measure the relative position between the points $x$ and $y$. And the angle $\theta(x,y)$ that appears in the coefficient $\tau^n_s(\theta(x,y))$ is the amount of deviation from the horizontal state between points $x$ and $y$ in the Heisenberg group.

The paper is organized as follows. In Section 2 we recall some basic facts on left-invariant vector fields, Carnot-Carath\'{e}odory distance, geodesic and optimal transport in $\mathbb{H}^n$. In Section 3, we introduce the curvature-dimension condition, the Brunn-Minkowski inequality and the strong Brunn-Minkowski inequality in the Heisenberg group, and discuss the interaction between the Brunn-Minkowski inequality and the strong Brunn-Minkowski inequality. The proof of the strong Brunn-Minkowski inequality in the Heisenberg group is also given in this section. We will give our main result, Theorem 4.1, and its proof in the last Section.

\section{\bf Preliminaries}

\noindent{\bf 2.1 The Heisenberg group}

We are now ready to introduce the following significant notations and definitions in the Heisenberg group $\mathbb{H}^n$. The Lie algebra of the Heisenberg group is spanned by the following left-invariant vector fields:
$$X_j=\partial_{\xi_j}+2\eta_j\partial_t,~~Y_j=\partial_{\eta_j}-2\xi_j\partial_t,~j\in\{1,\ldots,n\},~~T=\partial_t.$$
The horizontal distribution at a point $x\in\mathbb{H}^n$ is defined by
$$H_x=span\{X_i(x),Y_i(x): i=1,\ldots,n\}.$$
The horizontal distribution is nonintegrable. In fact there holds $[X_i,Y_i]=-4T\neq0$ for any $i=1,\ldots,n$, and all other commutators vanish. Roughly speaking, the Carnot-Carath\'{e}odory distance between two points $x$ and $y$ is the infimum of the lengths of the horizontal curves connecting $x$ and $y$. A horizontal curve whose derivative $\gamma'(s)$ is spanned by
$$\{X_1(\gamma(s)),\ldots,X_n(\gamma(s)),Y_1(\gamma(s)),\ldots,Y_n(\gamma(s))\},$$
in almost every point $s\in[0,1]$. The horizontal length of this curve is then
$$Length_{cc}(\gamma)=\int_0^1\|\gamma'(s)\|_\mathbb{H}ds,$$
where $\Big\|\sum \limits^n \limits_{k=1}(a_kX_k+b_kY_k)\Big\|_\mathbb{H}^2=\sum \limits^n \limits_{k=1}(a_k^2+b_k^2)\leq1$. Then the Carnot-Carath\'{e}odory distance $d(x,y)$ between $x$ and $y$ of $\mathbb{H}^n$ is
$$d(x,y):=inf \int_0^1\|\gamma'(s)\|_\mathbb{H}ds=inf Length_{cc}(\gamma),$$
where the infimum is taken over all horizontal curve connecting $x$ and $y$, and the length minimizing horizontal curves joining pairs of points in the Heisenberg group is called geodesic. Besides, the Carnot-Carath\'{e}odory unit ball centered at point $x\in\mathbb{H}^n$ is defined by $\mathbf{B}(x)=\{y\in\mathbb{H}^n: d(x,y)\leq1\}$.

The equations of the local geodesics of $\mathbb{H}^n$ have been known since Gaveau's paper \cite {G0}, and Ambrosio and Rigot \cite {AR} also study the globality of local geodesics. Because the Carnot-Carath\'{e}odory distance, hence the geodesics are left invariant, any geodesic can be left translated back to $0_{\mathbb{H}^n}$, it is enough to know the equations of the geodesics passing through $0_{\mathbb{H}^n}$. For every $(\chi,\theta)\in\mathbb{C}^n\times\mathbb{R}$ we mean the curve $\gamma_{\chi,\theta}: [0,1]\rightarrow\mathbb{H}^n$ defined by
\begin{align*}
\gamma_{\chi,\theta}(s)=
\begin{cases}
 \Big(i\frac{e^{i\theta s-1}}{\theta}\chi,2|\chi|^2\frac{\theta s-\sin(\theta s)}{\theta^2}\Big)\ \ \ \ \ \ \ \ \ if~ \theta\neq0;\\
(s\chi,0) \ \ \ \ \ \ \ \ \ \ \ \ \ \ \ \ \ \ \ \ \ \ \ \ \ \ \ \ \ \ \ \ \ if ~\theta=0.
\end{cases}
\end{align*}
Here we get the projected curve $\chi$ and the angle $\theta$ from the Dido's question, $|\chi|$ is $\sqrt{|\chi_1|^2+\ldots+|\chi_n|^2}$ and $(\chi,\theta)\in(\mathbb{C}^n\setminus \{0_{\mathbb{C}^n}\})\times[-2\pi,2\pi]$, the paths $\gamma_{\chi,\theta}$ are length-minimizing geodesics in $\mathbb{H}^n$ joining~~$0_{\mathbb{H}^n}$ and $\gamma_{\chi,\theta}(1)$. In this paper, we let
$$\Gamma_s(\chi,\theta):=\gamma_{\chi,\theta}(s).\eqno(2.1)$$
In particular, we will make use of it for $s=1$. We note that $\Gamma_s(\chi,\theta)$ is $\Gamma_1(s\chi,s\theta)$. If $\theta\in(-2\pi,2\pi)$ the geodesics connecting $0_{\mathbb{H}^n}$ and $\gamma_{\chi,\theta}(1)\neq0_{\mathbb{H}^n}$ are unique, while for $\theta\in\{-2\pi,2\pi\}$ the uniqueness fails. Obviously, the inverse map  $\Gamma_1^{-1}(\gamma_{\chi,\theta}(1))=(\chi,\theta)\in\mathbb{C}^n\times[-2\pi,2\pi]$.

For $s\in[0,1]$, the Heisenberg distortion coefficients $\tau^n_s: [0,2\pi]\rightarrow[0,+\infty]$ defined by
\begin{align*}
\tau^n_s(\theta)=
\begin{cases}
+\infty  \ \ \ \ \ \ \ \ \ \ \ \ \ \ \ \ \ \ \ \ \ \ \ \ \ \ \ \ \ \ \ \ \ \ \ \ \ \ \ \ \ \ \ \ \ \ if~~~ \theta=2\pi;\\
s^\frac{1}{2n+1}\Big(\frac{\sin\frac{\theta s}{2}}{\sin\frac{\theta}{2}} \Big)^\frac{2n-1}{2n+1} \Big(\frac{\sin\frac{\theta s}{2}-\frac{\theta s}{2}\cos \frac{\theta s}{2}}{\sin\frac{\theta}{2}-\frac{\theta}{2}\cos\frac{\theta}{2}}\Big)^\frac{1}{2n+1}~~if~~~\theta\in(0,2\pi); \label{eqsystem1}\\
s^\frac{2n+3}{2n+1} \ \ \ \ \ \ \ \ \ \ \ \ \ \ \ \ \ \ \ \ \ \ \ \ \ \ \ \ \ \ \ \ \ \ \ \ \ \ \ \ \ \ \ \ \  if ~~~\theta=0.
\end{cases}
\end{align*}
The function $\theta\mapsto\tau^n_s(\theta)$ is increasing on $[0,2\pi]$. Specifically for $\theta\rightarrow 2\pi$,
$\tau^n_s(\theta)\rightarrow+\infty$, and
$$\tau^n_s(\theta)\geq\tau^n_s(0)=s^\frac{2n+3}{2n+1}~~~for~every~\theta\in[0,2\pi],~s\in[0,1].$$
Here, if $x,y\in\mathbb{H}^n$, $x\neq y$ we set $\theta(x,y)=\theta(y,x)$. If $x=y$, $\theta(x,y)=0$.

\noindent{\bf 2.2 Optimal transport}

The problem of optimal transportation raised by G. Monge, it has received a lot of attention in recent years due to its increasing application (see for instance \cite {E,RR,V}). In the Monge problem, given two Borel probability measures $\mu$, $\nu$ in a topological space $X$ and minimize
$$\Psi\mapsto\int_Xc(x,\Psi(x))d\mu(x)$$
among all transport maps $\Psi$, i.e. all Borel maps $\Psi: X\rightarrow X$ such that $\Psi_\sharp\mu=\nu$ ($\Psi$ pushes $\mu$ into $\nu$, $\nu(B)=\mu(\Psi^{-1}(B))$ for all Borel set $B$), $c(x,y)$ is represents the cost for shipping a unit of mass from $x$ to $y$. In 1942, Kantorovich \cite {K0} extended the transport mapping into a transport planning: consider the
larger class of transport plannings, probability measures $\pi$ in $X\times X$ such that $p_{0\sharp}\pi=\mu$ and $p_{1\sharp}\pi=\nu$, here $p_i: X\times X\rightarrow X$ denotes the projection on the $i$-th factor ($i=0,1$), and to minimize
$$\pi\mapsto\int_{X\times X}c(x,y)d\pi(x,y)$$
in the class of transport plannings. Note that any transport map $\Psi$ induces via the map $(id\times\Psi)(x):=(x,\Psi(x))$ a planning $\pi$ with the same cost, defined by $(id\times\Psi)_\sharp\mu$. Over the years, optimal transmission has been extended to the Heisenberg group, next we give some representation about optimal transport on Heisenberg group, for more details, see \cite {AR,FJ,J1}.

A metric measure space is a triple $(\mathbb{H}^n,d,\mathcal L^{2n+1})$, where $d$ is the Carnot-Carath\'{e}odory distance between two points, $\mathcal L^{2n+1}$ is a Haar measure on $\mathbb{H}^n$ and $(\mathbb{H}^n,d)$ is a Polish metric space (complete and separable). We denote by $\mathcal{M}_+(\mathbb{H}^n)$ the set of all positive and finite Borel measures on metric space $(\mathbb{H}^n,d)$. Let $\mathcal{P}(\mathbb{H}^n)(\subset\mathcal{M}_+(\mathbb{H}^n))$ denote the set of all positive and finite probability measures on $(\mathbb{H}^n,d)$, i.e. if $\mu$ is a  Borel probability measures on $\mathbb{H}^n$, $\mu(\mathbb{H}^n)=1$. $\mathcal{P}_2(\mathbb{H}^n)$ denote the space of probability measures in $\mathbb{H}^n$ with finite second-order moments, endowed with the Wasserstein metric between $\mu,\nu\in\mathcal{P}_2(\mathbb{H}^n)$
$$W_2(\mu,\nu):=\Big(\min \limits_{\Upsilon}\int_{\mathbb{H}^n\times\mathbb{H}^n }d(x,y)^2d\Upsilon(x,y)\Big)^{1/2}, \eqno(2.2)$$
here probability measures $\Upsilon$ varies in the family $\Sigma(\mu,\nu)$ of admissible transport plannings, and $\Upsilon$ with fixed marginals $p_{0\sharp}\Upsilon=\mu$ and $p_{1\sharp}\Upsilon=\nu$ ($p_0,~p_1: \mathbb{H}^n\times\mathbb{H}^n\rightarrow\mathbb{H}^n$ being the canonical projections on the first and on the second factor) can be written as a coupling induced by a transport map, i.e. whether there is a measurable map $\Pi: \mathbb{H}^n\rightarrow\mathbb{H}^n$ such that $\Upsilon=(id\otimes \Pi)_\sharp\mu$.
 We set the family of measures
$$\mu_s:=(\Pi_s)_\sharp\mu,~s\in[0,1]~~with~~\Pi_{s}(x):=x\ast exp_\mathbb{H}(-sX\varphi-isY\varphi-s^2Z\varphi),$$
is a geodesic in $\mathcal{P}_2(\mathbb{H}^n)$ connecting $\mu$ and $\nu$. Here $ exp_\mathbb{H}$ is the sub-Riemannian exponential, $\varphi$ is a $c$-concave map. Since both $\mu$ and $\nu$ are absolutely continuous, there exists also an optimal transport map $\Phi$ from $\nu$ to $\mu$, and it is well-known that $\Phi$ is an inverse for $\Pi$ a.e., that is 
$$\Phi\circ\Pi=id~~\mu-a.e.,~~~~\Pi\circ\Phi=id~~\nu-a.e.$$
When the probability measures in $\mathcal{P}_2(\mathbb{H}^n)$ are absolutely continuous with respect to $\mathcal L^{2n+1}$ denoted by $\mathcal{P}^{ac}(\mathbb{H}^n,\mathcal L^{2n+1})$. The minimum in (2.2) is always attained, the admissible plans realizing it are called optimal transport plans and the set that contains all of them is denoted by $Opt(\mu,\nu)$.

We known $W_2$ is a complete and separable distance on $\mathcal{P}_2(\mathbb{H}^n)$. The convergence of Wasserstein distance is expressed as follows:
$$\mu_n \stackrel{^{W_2}}\rightarrow \mu~~~\Leftrightarrow~~~\mu_n\rightharpoonup\mu~and~\int d(x_0,x)^2d\mu_n~\rightarrow~\int d(x_0,x)^2d\mu~~~\forall~x_0\in \mathbb{H}^n,$$
here $\rightharpoonup$ denotes the weak convergence of measures.

Let $C([0,1],\mathbb{H}^n)$ be the set of continuous functions from $[0,1]$ to $\mathbb{H}^n$, and for $\gamma\in C([0,1],\mathbb{H}^n)$, $s\in[0,1]$ define the $s$-evaluation map as $T_s:C([0,1],\mathbb{H}^n)\rightarrow\mathbb{H}^n$, $T_s(\gamma):=\gamma(s)$. Let Geo($\mathbb{H}^n$) denote the space of constant speed geodesics in $(\mathbb{H}^n,d)$ parametrized on $[0,1]$. $(\mathbb{H}^n,d)$ is said to be geodesic if each pair of points can be connected by a geodesic. If~$(\mathbb{H}^n,d)$ is geodesic metric space then $(\mathcal{P}_2(\mathbb{H}^n),W_2)$ is also geodesic metric space. In the space of probability measures $\mathcal{P}(\mathbb{H}^n)$, every measure $\eta\in\mathcal{P}(C([0,1],\mathbb{H}^n))$ induces the curve $[0,1]\ni s\mapsto\mu_s=(T_s)_\sharp\eta$. More accurately, given a pair of measures $\mu$, $\nu\in\mathcal{P}_2(\mathbb{H}^n)$, the curve $\{\mu_s\}_{s\in [0,1]}\subset\mathcal{P}_2(\mathbb{H}^n)$ connecting them is a Wasserstein geodesic if and only if there exists $\eta\in\mathcal{P}(C([0,1],\mathbb{H}^n))$ inducing curve $\{\mu_s\}_{s\in [0,1]}$, which is concentrated on  Geo($\mathbb{H}^n$) and satisfies $(T_0,T_1)_\sharp\eta\in Opt(\mu,\nu)$. In this case, $\eta$ is called the optimal geodesic plan between $\mu$ and $\nu$ and this will be denoted as $\eta\in OptGeo(\mu,\nu)$.

\noindent{\bf Definition 2.1}(\cite {K}).~~{\it A set of geodesics $G\subset Geo(\mathbb{H}^n)$ is called non-branching if for all geodesics $\gamma$, $\eta\in G$ with $restr_0^s \gamma= restr_0^s \eta$ for some $s\in(0,1)$ it holds $\gamma\equiv\eta$.}

A metric measure space $(\mathbb{H}^n,d,\mathcal L^{2n+1})$ is said to be {\it essentially non-branching} if for every pair
absolutely continuous probability measures $\mu_0$, $\mu_1\in\mathcal{P}^{ac}(\mathbb{H}^n,\mathcal L^{2n+1})$, every optimal geodesic plan $\eta$ connecting them is concentrated on the non-branching set of geodesics.

\noindent{\bf Definition 2.2}(\cite {K}).~~{\it The measure $\mathcal L^{2n+1}$ is said to be qualitatively non-degenerate if for all $R>0$ and $x_0\in\mathbb{H}^n$ there is a function $f_{R,x_0}: (0,1)\rightarrow(0,\infty)$ with
$$\mathop{\lim \sup}\limits_{s\rightarrow 0}f_{R,x_0}(s)>\frac{1}{2},$$
such that for all measurable set $A\subset \mathbf{B}_R(x_0)$ and $x\in \mathbf{B}_R(x_0)$ and $s\in(0,1)$ it holds
$$\mathcal L^{2n+1}(Z_s(A,x))\geq f_{R,x_0}(s)\cdot\mathcal L^{2n+1}(A),$$
here  the s-intermediate set $Z_s(A,x)$ is defined in section 3.}

\noindent{\bf Definition 2.3}(\cite {K}).~~{\it A metric measure space $(\mathbb{H}^n,d,\mathcal L^{2n+1})$ has good transport behavior if for all $\mu$, $\nu\in\mathcal{P}_2(\mathbb{H}^n)$ with $\mu\ll\mathcal L^{2n+1}$ any optimal transport plan between $\mu$ and $\nu$ is induced by a map. $(\mathbb{H}^n,d,\mathcal L^{2n+1})$ has strong interpolation property if for all $\mu,\nu\in\mathcal{P}_2(\mathbb{H}^n)$ with $\mu\ll\mathcal L^{2n+1}$, there exists a unique optimal geodesic plan $\eta\in OptGeo(\mu,\nu)$ and is induced by a map and such that $(T_s)_\sharp\eta\ll\mathcal L^{2n+1}$ for $s\in(0,1)$.}

\noindent{\bf Proposition 2.4}(\cite {K}).~~{\it If~$(\mathbb{H}^n,d,\mathcal L^{2n+1})$ is a proper, geodesic, essentially non-branching metric measure space and $\mathcal L^{2n+1}$ is qualitatively non-degenerate, then $(\mathbb{H}^n,d,\mathcal L^{2n+1})$ has both the good transport behaviour and the strong interpolation property.}

\section{\bf Some Results on the Heisenberg group}

$\quad$~~Let us recall some facts below. A function $U:[0,\infty)\rightarrow\mathbb{R}$ defines the R\'{e}nyi entropy (when $U(t)=-t^{1-\frac{1}{2n+1}}$) of an absolutely continuous measure $\mu$ with respect to~~$\mathcal L^{2n+1}$ on $\mathbb{H}^n$ as (see \cite {BKS})
$$Ent_U(\mu|\mathcal L^{2n+1})=\int_{\mathbb{H}^n}U(\rho(x))d\mathcal L^{2n+1}(x)=\int_{\mathbb{H}^n}-\rho(x)^{1-\frac{1}{2n+1}}d\mathcal L^{2n+1}(x),$$
where $\rho=\frac{d\mu}{d\mathcal L^{2n+1}}$ is the density function of measure $\mu$. There are many types of entropy, as can be seen in \cite {G1,W}.

\noindent{\bf Theorem 3.1} (CD condition on $\mathbb{H}^n$ \cite {BKS}).~~{\it Let $s\in[0,1]$ and assume that $\mu_0$ and $\mu_1$ are two compactly supported, Borel probability measures, both absolutely continuous with respect to $\mathcal L^{2n+1}$ on $\mathbb{H}^n$ with densities $\rho_0$ and $\rho_1$, respectively, i.e. $\mu_0=\rho_0\mathcal L^{2n+1}$, $\mu_1=\rho_1\mathcal L^{2n+1}\in\mathcal{P}^{ac}(\mathbb{H}^n,\mathcal L^{2n+1})$. If exists a $W_2$-geodesic $\eta\in\mathcal{P}(Geo(\mathbb{H}^n))$ connecting $\mu_0$ and $\mu_1$, such that $\mu_s=(T_s)_\sharp\eta=\rho_s\mathcal L^{2n+1}\ll\mathcal L^{2n+1}$, the following entropy inequality holds:
$$Ent_{2n+1}(\mu_s|\mathcal L^{2n+1})\leq-\int_{\mathbb{H}^n\times\mathbb{H}^n}\Big[\tau^n_{1-s}(\theta_x)\rho_0(x)^{-\frac{1}{2n+1}}+\tau^n_s(\theta_x)\rho_1(y)^{-\frac{1}{2n+1}}
\Big]d\Upsilon(x,y),\eqno(3.1)$$
here $\theta_x=\theta_{x,T(x)}=|\theta(x,T(x))|$ and $\Upsilon=(T_0,T_1)_\sharp\eta\in Opt(\mu_0,\mu_1)$.}

For convenience, we sometimes use $\mathcal{F}^n_s(\Upsilon|\mathcal L^{2n+1})$ to represent the right side of (3.1), this means that
$$\mathcal{F}^n_s(\Upsilon|\mathcal L^{2n+1})=-\int_{\mathbb{H}^n\times\mathbb{H}^n}\Big[\tau^n_{1-s}(\theta_x)\rho_0(x)^{-\frac{1}{2n+1}}+\tau^n_s(\theta_x)\rho_1(y)^{-\frac{1}{2n+1}}
\Big]d\Upsilon(x,y).$$

\noindent{\bf Definition 3.2} (Bounded probability measure).~~{\it A probability measure $\mu\in\mathcal{P}^{ac}(\mathbb{H}^n,\mathcal L^{2n+1})$ is said to be bounded if it has bounded support and density bounded from above and below away from zero. A subset $A\subset\mathcal{P}^{ac}(\mathbb{H}^n,\mathcal L^{2n+1})$ is uniformly bounded if there exist a bounded set $K$ and two constants $D>d>0$ such that for every $\mu=\rho\mathcal L^{2n+1}\in A$, $spt(\mu)=K$ and $d\leq\rho\leq D$ $\mathcal L^{2n+1}$-almost everywhere on $K$.}

Next we introduce the concept of the step measure. The step measure is one of the simpler measures in structure, any bounded measure can be approximated by a suitable step measure, so we proceed to complete the proof by applying approximation.

\noindent{\bf Definition 3.3}.~~{\it The measure $\mu\in\mathcal{P}_2(\mathbb{H}^n)$ is called a step measure if $\mu$ can be written as the finite sum of measures with constant density on $\mathcal L^{2n+1}$, namely
$$\mu=\sum \limits^N \limits_{i=1} \lambda_i \mathcal L^{2n+1}_{A_i},$$
for each $i$, $\lambda_i\in\mathbb{R}$ and mutually disjoint open Borel set $\{A_i\}_{i=1,\ldots,N}$ with $0<\mathcal L^{2n+1}(A_i)<\infty$.}

\noindent{\bf Lemma 3.4}.~~{\it Let $\mu=\rho\mathcal L^{2n+1}\in\mathcal{P}^{ac}(\mathbb{H}^n,\mathcal L^{2n+1})$ be bounded, then there exists a sequence of step measures $\{\mu_m=\rho_m\mathcal L^{2n+1}\}_{m\in \mathbb{N}}$ $W_2$-convergent to $\mu$, such that $\{\mu_m\}_{m\in \mathbb{N}}\bigcup\{\mu\}$ is uniformly bounded and $\rho_m^{-\frac{1}{2n+1}}\rightarrow\rho^{-\frac{1}{2n+1}}$ in $L^1$.}

\noindent{\bf Lemma 3.5}.~~{\it Let $(\mathbb{H}^n,d,\mathcal L^{2n+1})$  be the metric measure space, two bounded measure $\mu, \nu\in\mathcal{P}^{ac}(\mathbb{H}^n,\mathcal L^{2n+1})$ and $\Upsilon$ is the unique optimal transport plan between them. Assume $\{\mu_m\}_{m\in \mathbb{N}}$, $\{\nu_m\}_{m\in \mathbb{N}}\in\mathcal{P}^{ac}(\mathbb{H}^n,\mathcal L^{2n+1})$ be the approximating sequences provided by Lemma 3.4, letting $\Upsilon_m$ is the unique optimal transport plan between $\mu_m$ and $\nu_m$ and weakly converging to $\Upsilon$, then
$$\mathop{\lim \sup}\limits_{m\rightarrow \infty}\mathcal{F}^n_s(\Upsilon_m|\mathcal L^{2n+1})\leq\mathcal{F}^n_s(\Upsilon|\mathcal L^{2n+1}),$$
for all $s\in[0,1]$.}

\noindent{\bf Proposition 3.6}.~~{\it Assume $(\mathbb{H}^n,d,\mathcal L^{2n+1})$ is an essentially non-branching metric measure space and $\mathcal L^{2n+1}$ is qualitatively non-degenerate. Then, $(\mathbb{H}^n,d,$
$\mathcal L^{2n+1})$ satisfies the CD condition if the inequality (3.1) holds for any pair of bounded probability measures.}

The proof of Lemma 3.4, Lemma 3.5 and Proposition 3.6 can be found in Ref. \cite{MPR}. The case in Riemannian spaces was proved in \cite{MPR}, and only the distortion coefficients are different in the Heisenberg group, and it is easy to prove that the above conclusions are appropriate in a similar way. 

Notation: Given a Borel set $A\subset\mathbb{H}^n$ and $0<\mathcal L^{2n+1}(A)<\infty$, we will denote by $\mathcal L^{2n+1}_A$ the formalization  of the reference measure to the set $A$, that is
$$\mathcal L^{2n+1}_A=\frac{\mathcal L^{2n+1}|A}{\mathcal L^{2n+1}(A)}.$$

For two nonempty measurable sets $A,B\subset\mathbb{H}^n$, let $\Theta_{A,B}$ represents a typical Heisenberg quantity indicating a lower bound of the deviation of an essentially horizontal position of the sets A and B:
$$\Theta_{A,B}=\sup \limits_{A_0,B_0}\inf\limits_{(x,y)\in A_0\times B_0}\{|\theta|\in[0,2\pi]: (\chi,\theta)\in\Gamma_1^{-1}(x^{-1}\ast y)\},\eqno(3.2)$$
here $A_0$ and $B_0$ are nonempty, full measure subsets of $A$ and $B$, respectively (see \cite{BKS}). According to Figalli and Juillet \cite{FJ}, for $s\in[0,1]$, $Z_s(A,B)$ denotes the s-intermediate set associated to the nonempty measurable sets $A,B\subset\mathbb{H}^n$ with respect to the Carnot-Carath\'{e}odory metric $d$. Note that $(\mathbb{H}^n,d)$ is a geodesic metric space, thus $Z_s(x,y)\neq\emptyset$ for each pair of points $x,y\in\mathbb{H}^n$. According to Ambrosio and Rigot \cite {AR}, there exists a unique optimal transport map $\Pi:\mathbb{H}^n\rightarrow\mathbb{H}^n$ transporting $\mu_0$ to $\mu_1$ associated to the cost function $d^2$, if $\Pi_s$ denotes the interpolant optimal transport map, then
$$\Pi_s(x)=Z_s(x,\Pi(x)),~Z_s(A,B)=\bigcup_{(x,y)\in A\times B}Z_s(x,y)~~~~x,~y\in\mathbb{H}^n.\eqno(3.3)$$
Generally, for $x\in\mathbb{H}^n$ and $A\subset\mathbb{H}^n$, we record as: $Z_s(x,A):=Z_s(\{x\},A)$, $Z_s(A,x):=Z_s(A,\{x\})$.

\noindent{\bf Theorem 3.7} (Brunn-Minkowski inequality on $\mathbb{H}^n$ \cite {BKS}).~~{\it Let $(\mathbb{H}^n,d,\mathcal L^{2n+1})$ be a metric measure space, for every two nonempty Borel sets $A,B \subset spt(\mathcal L^{2n+1})$ and $s\in[0,1]$, the following geodesic Brunn-Minkowski inequality BMI holds:
$$\mathcal L^{2n+1}(Z_s(A,B))^\frac{1}{2n+1}\geq\tau^n_{1-s}(\Theta_{A,B})\mathcal L^{2n+1}(A)^\frac{1}{2n+1}+\tau^n_s(\Theta_{A,B})\mathcal L^{2n+1}(B)^\frac{1}{2n+1}, \eqno(3.4)$$
here $\Theta_{A,B}$ is defined by (3.2).}

We consider the outer measure when $Z_s(A,B)$ is not measurable. When $A^{-1}\ast B\subset C=\{0_\mathbb{C}^n\}\times\mathbb{R}$, $\Theta_{A,B}=2\pi$ and $\mathcal L^{2n+1}(A)=\mathcal L^{2n+1}(B)=0$.

\noindent{\bf Lemma 3.8}.~~{\it Assume that the metric measure space $(\mathbb{H}^n,d,\mathcal L^{2n+1})$ satisfies BMI, then $(spt(\mathcal L^{2n+1}),d)$ is a Polish, geodesic and proper metric space, and also $\mathcal L^{2n+1}$ is a Radon measure.}

{\it Proof.}~~Since $(\mathbb{H}^n,d)$ is Polish and $spt(\mathcal L^{2n+1})$ is closed, so the metric space $(spt(\mathcal L^{2n+1}),d)$ is Polish as well. It satisfies BMI to know that any bounded set is totally bounded, therefore it is proper and $\mathcal L^{2n+1}$ is Radon, so that being a locally finite measure on a locally compact space. Next, We prove $(spt(\mathcal L^{2n+1}),d)$ is a length space: if $x,y\in spt(\mathcal L^{2n+1})$ and $\varepsilon>0$, let $A_\varepsilon: =\mathbf{B}_\varepsilon(x)\bigcap spt(\mathcal L^{2n+1})$ and $B_\varepsilon: =\mathbf{B}_\varepsilon(y)\bigcap spt(\mathcal L^{2n+1})$. Applying BMI we deduce that $\mathcal L^{2n+1}(Z_\frac{1}{2}(A_\varepsilon,B_\varepsilon))>0$, therefore there exists $z\in Z_\frac{1}{2}(A_\varepsilon,B_\varepsilon)\bigcap spt(\mathcal L^{2n+1})$ such that
$$d(x,z);d(z,y)\leq \frac{1}{2}d(x,y)+\varepsilon,$$
so $(spt(\mathcal L^{2n+1}),d)$ is a length space (see \cite{S1}). Finally, a complete, proper and length space is geodesic.
\hfill ${\square}$

Combining Lemma 3.8 and Proposition 2.4, we get that the metric measure space $(\mathbb{H}^n,d,\mathcal L^{2n+1})$ satisfies Brunn-Minkowski inequality has the following property:

\noindent{\bf Corollary 3.9.}~~{\it Assume that $(\mathbb{H}^n,d,\mathcal L^{2n+1})$ is an essentially non-branching metric measure space supporting the Brunn-Minkowski inequality, then $(\mathbb{H}^n,d,\mathcal L^{2n+1})$ has the good
transport behavior and the strong interpolation property.}

In this paper, we study a stronger version of the Brunn-Minkowski inequality on $\mathbb{H}^n$, which is involves the optimal transport interpolation.

\noindent{\bf Theorem 3.10} (Strong Brunn-Minkowski inequality on $\mathbb{H}^n$).~~{\it Let $(\mathbb{H}^n,d,\mathcal L^{2n+1})$ be an essentially non-branching metric measure space and $\mathcal L^{2n+1}$ is qualitatively non-degenerate. If for every two nonempty Borel sets $A,B \subset spt(\mathcal L^{2n+1})$ and $0<\mathcal L^{2n+1}(A)$, $\mathcal L^{2n+1}(B)<\infty$, exists $\eta\in OptGeo(\mathcal L^{2n+1}_A,\mathcal L^{2n+1}_B)$, such that the following inequality holds for $s\in[0,1]$:
$$\mathcal L^{2n+1}(spt((T_s)_\sharp\eta))^\frac{1}{2n+1}\geq\tau^n_{1-s}(\Theta_{A,B})\mathcal L^{2n+1}(A)^\frac{1}{2n+1}+\tau^n_s(\Theta_{A,B})\mathcal L^{2n+1}(B)^\frac{1}{2n+1}. \eqno(3.5)$$
We call inequality (3.5) the SBMI and $\Theta_{A,B}$ is introduced in (3.2).}

\noindent{\bf Lemma 3.11} (Borell-Brascamp-Lieb inequality on $\mathbb{H}^n$ \cite {BKS}).~~{\it For any $s\in[0,1]$ and $p\geq-\frac{1}{2n+1}$, let $f,g,h: \mathbb{H}^n\rightarrow[0,\infty)$ be integrable functions with the property that for all $(x,y)\in\mathbb{H}^n\times\mathbb{H}^n, z\in Z_s(x,y)$,
$$h(z)\geq M_s^p(\frac{f(x)}{(\widetilde{\tau}^n_{1-s}(\theta(y,x)))^{2n+1}},\frac{g(y)}{(\widetilde{\tau}^n_{s}(\theta(x,y)))^{2n+1}}),$$
then the following inequality holds:
$$\int_{\mathbb{H}^n}h\geq M_s^{\frac{p}{1+(2n+1)p}}(\int_{\mathbb{H}^n}f,\int_{\mathbb{H}^n}g).$$}
Here $\widetilde{\tau}_s^n(\theta)=s^{-1}\tau_s^n(\theta)$, $p\in\mathbb{R}\cup\pm\infty$ and $a,b\geq0$, we consider the p-mean
\begin{align*}
M_s^p(a,b)=
\begin{cases}
((1-s)a^p+sb^p)^{\frac{1}{p}}  \ \ \ \ \ \ \ \ \ \ \ \ \ \ \ \ \ \ \ \ \ \ if~~~ ab\neq0;\\
0 \ \ \ \ \ \ \ \ \ \ \ \ \ \ \ \ \ \ \ \ \ \ \ \ \ \ \ \ \ \ \ \ \ \ \ \ \ \ \ \ \ \ \ \ if ~~~ab=0.
\end{cases}
\end{align*}
{\it Proof of Theorem 3.10.}~~As $(\mathbb{H}^n, d, \mathcal L^{2n+1})$ is an non-branching metric measure space and $\mathcal L^{2n+1}$ is qualitatively non-degenerate, then $(\mathbb{H}^n, d, \mathcal L^{2n+1})$ has the good transport behavior and the strong interpolation property. For any two Borel sets $A,B \subset spt(\mathcal L^{2n+1})$, $0<\mathcal L^{2n+1}(A)$, $\mathcal L^{2n+1}(B)<\infty$, by optimal transport we know that $(T_s)_\sharp\eta$ is the measure corresponding to $Z_s(A,B)$ and $spt((T_s)_\sharp\eta)\subset Z_s(A,B)$.

We recall
$$\Theta_{A,B}=\sup \limits_{A_0,B_0}\inf\limits_{(x,y)\in A_0\times B_0}\{|\theta|\in[0,2\pi]: (\chi,\theta)\in\Gamma_1^{-1}(x^{-1}\ast y)\},$$
here $A_0$ and $B_0$ are nonempty, full measure subsets of $A$ and $B$, respectively. If $\Theta_{A,B}=2\pi$, for any $(x,y)\in A\times B$ we have $x^{-1}\ast y\in\Gamma_1(\chi,\pm 2\pi)\subset C={(\zeta,t)\in\mathbb{H}^n:\zeta=0}$, $A^{-1}\ast B\subset C$, then $\mathcal L^{2n+1}(A^{-1}\ast B)=0$. From the multiplicative Brunn-Minkowski inequality
$$\mathcal L^{2n+1}(A^{-1}\ast B)^{\frac{1}{2n+1}}\geq\mathcal L^{2n+1}(A^{-1})^{\frac{1}{2n+1}}+\mathcal L^{2n+1}(B)^{\frac{1}{2n+1}},$$
we get $\mathcal L^{2n+1}(A)=\mathcal L^{2n+1}(B)=0$, therefore $\Theta_{A,B}$ is not eequal to $2\pi$, i.e. $\Theta_{A,B}<2\pi$.

For $\Theta_{A,B}<2\pi$ and $s\in[0,1]$, let
$$c_1^s=\sup \limits_{A_0,B_0}\inf\limits_{(x,y)\in A_0\times B_0}\widetilde{\tau}^n_{1-s}(\theta(y,x))=\widetilde{\tau}^n_{1-s}(\theta(A,B)),$$
$$c_2^s=\sup \limits_{A_0,B_0}\inf\limits_{(x,y)\in A_0\times B_0}\widetilde{\tau}^n_{s}(\theta(x,y))=\widetilde{\tau}^n_{s}(\theta(A,B)).$$
Where $A_0$ and $B_0$ are nonempty, full measure subsets of $A$ and $B$, respectively. $\widetilde{\tau}_s^n(\theta)=s^{-1}\tau_s^n(\theta)$, $\tau_s^n(\theta)$ is increasing on $[0,2\pi)$ and $0<c_1^s,c_2^s<+\infty$. Let $\mathcal L^{2n+1}(A)\neq 0\neq\mathcal L^{2n+1}(B)$, $p=+\infty$, $f(x)=(c_1^s)^{2n+1}\mathds{1}_A(x)$, $g(y)=(c_2^s)^{2n+1}\mathds{1}_B(y)$ and $h(z)=\mathds{1}_{spt((T_s)_\sharp\eta)}(z)$, use Lemma 3.11, we have
$$\mathcal L^{2n+1}(spt((T_s)_\sharp\eta))\geq M_s^{\frac{1}{2n+1}}((c_1^s)^{2n+1}\mathcal L^{2n+1}(A),(c_2^s)^{2n+1}\mathcal L^{2n+1}(B))\ \ \ \ \ \ \ \ \ \ \ \ \ \ \ \ \ \ \ \ \ \ \ \ \ \ \ \ \ \ \ \ \ \ \ $$
$$\ \ \ \ \ \ \ \ \ \ \ \ =(\tau^n_{1-s}(\Theta_{A,B})\mathcal L^{2n+1}(A)^\frac{1}{2n+1}+\tau^n_s(\Theta_{A,B})\mathcal L^{2n+1}(B)^\frac{1}{2n+1})^{2n+1},$$
which concludes the proof. \hfill ${\square}$

\noindent{\bf Remark:}~(i) Since $spt((T_s)_\sharp\eta)\subset Z_s(A,B)$, when the SBMI in metric measure space $(\mathbb{H}^n,d$, $\mathcal L^{2n+1})$ holds, the BMI also holds. By the Corollary 3.9 and Theorem 3.10, it follows that an essentially non-branching metric measure space $(\mathbb{H}^n,d,\mathcal L^{2n+1})$ supporting SBMI has the good transport behavior and the strong interpolation property. Thus, given two nonempty Borel sets $A,B\subset spt(\mathcal L^{2n+1})$ with finite and positive measure, then there exists a unique optimal geodesic plan $\eta\in OptGeo(\mathcal L^{2n+1}_A$, $\mathcal L^{2n+1}_B)$ depending only on the sets $A$ and $B$.

~(ii) Let $\lambda>0$, $\forall~x,~y\in\mathbb{H}^n$, since $(\delta_\lambda(x))^{-1}\ast\delta_\lambda(y)=\delta_\lambda(x^{-1}\ast y)$, for every Borel sets $A,B\in\mathbb{H}^n$ it turns out that $\Theta_{\delta_\lambda (A),\delta_\lambda (B)}=\Theta_{A,B}$. As a consequence, the strong
Brunn-Minkowski inequality and entropy inequality are invariant under the dilation of the sets.

\section{\bf Main Results and its Proof}

$\quad$~~In this section, we shall present our main result and its proof. The proof of necessity is given in Theorem 4.2, and the proof of sufficiency is given in Theorem 4.3.

\noindent{\bf Theorem 4.1}.~~{\it Let~$(\mathbb{H}^n, d, \mathcal L^{2n+1})$ be an essentially non-branching metric measure space and $\mathcal L^{2n+1}$ is qualitatively non-degenerate. Assume that two nonempty finite Borel sets $A,B \subset spt(\mathcal L^{2n+1})$, $\mu$ and $\nu$ are two bounded Borel probability measures attached to $A$ and $B$ respectively, and they are absolutely continuous with respect to $\mathcal L^{2n+1}$ on $\mathbb{H}^n$ with densities $\rho_0$ and $\rho_1$, then
$$(\mathbb{H}^n, d, \mathcal L^{2n+1})~supports~SBMI~if~and~only~if~it~satisfies~CD~condition.$$}

\noindent{\bf Theorem 4.2}.~~{\it Let the metric measure space $(\mathbb{H}^n,d,\mathcal L^{2n+1})$ be an essentially non-branching metric measure space and $\mathcal L^{2n+1}$ is qualitatively non-degenerate, assume that $(\mathbb{H}^n,d,\mathcal L^{2n+1})$ satisfies the curvature-dimension condition, then the strong Brunn-Minkowski inequality holds.}

{\it Proof.}~~Let us first give two Borel measurable sets $A,B \subset spt(\mathcal L^{2n+1})$ such that $0<\mathcal L^{2n+1}(A)$, $\mathcal L^{2n+1}(B)<\infty$, take the optimal geodesic plan $\eta\in OptGeo(\mathcal L^{2n+1}_A,\mathcal L^{2n+1}_B)$ satisfying inequality (3.1) and combined with the monotonicity of the Heisenberg distortion coefficients, we prove that
$$Ent_{2n+1}((T_s)_\sharp\eta)\leq-\int_{\mathbb{H}^n\times\mathbb{H}^n}\Big[\tau^n_{1-s}(\theta_x)\mathcal L^{2n+1}(A)^\frac{1}{2n+1}+\tau^n_s(\theta_x)\mathcal L^{2n+1}(B)^\frac{1}{2n+1}
\Big]d\Upsilon(x,y)$$
$$\ \ \ \ \ \ \ \leq-\Big[\tau^n_{1-s}(\theta_{A,B})\mathcal L^{2n+1}(A)^\frac{1}{2n+1}+\tau^n_s(\theta_{A,B})\mathcal L^{2n+1}(B)^\frac{1}{2n+1}
\Big].\eqno(4.1)$$
On the other side, recalling Jensen's inequality we get that
$$Ent_{2n+1}((T_s)_\sharp\eta)=-\int_{spt((T_s)_\sharp\eta)}\rho_s(x)^{1-\frac{1}{2n+1}}d\mathcal L^{2n+1}(x)\geq-\mathcal L^{2n+1}(spt((T_s)_\sharp\eta))^\frac{1}{2n+1},\eqno(4.2)$$
here $s\in[0,1]$, $\rho_s$ is the density of $(T_s)_\sharp\eta$ with respect to $\mathcal L^{2n+1}$. By inequalities (4.1) and (4.2), it is easy to obtain inequality (3.5), which concludes the proof.\hfill ${\square}$

A key idea for the proof of sufficiency is to prove the CD condition for a suitable subclass of bounded probability measures (the step measures). We known that $A,B \subset spt(\mathcal L^{2n+1})$ with finite and positive measure, the $SBMI$ would translate directly to an information on the entropy of the $s$-midpoint $\mu_s$ between $\mathcal L^{2n+1}_{A}$ and $\mathcal L^{2n+1}_{B}$, only if $\mu_s$ had constant density. At first, We will discuss in partition and establish the optimal transport coupling, the $s$-midpoint $\mu_s$ between $\mathcal L^{2n+1}_{A}$ and $\mathcal L^{2n+1}_{B}$ can be approximated in entropy of the step measure and has locally constant density. In fact, this partitioned discussion also applies when replacing the measures $\mathcal L^{2n+1}_{A}$ and $\mathcal L^{2n+1}_{B}$ with general step measures. In the end, we approximate the CD condition of all bounded measures by the CD condition for a class of step measures.

\noindent{\bf Theorem 4.3}.~~{\it Let $(\mathbb{H}^n,d,\mathcal L^{2n+1})$ be an essentially non-branching metric measure space and $\mathcal L^{2n+1}$ is qualitatively non-degenerate, for every pair of bounded measure $\mu,\nu\in\mathcal{P}^{ac}(\mathbb{H}^n$, $\mathcal L^{2n+1})$, if this space supporting the SBMI then the CD condition also holds.}

{\it Proof.}~~(i) We first demonstrate that the step measure $\mu_0,\mu_1\in\mathcal{P}^{ac}(\mathbb{H}^n,\mathcal L^{2n+1})$ with bounded support, if this space supporting the SBMI then the CD condition also holds.

From Corollary 3.9 and Theorem 3.10, we know that $(\mathbb{H}^n,d,\mathcal L^{2n+1})$ has the good transport behavior and the strong interpolation property. Therefore, letting $\mu_0$ and $\mu_1$ be two step measures, i.e.
$$\mu_0=\sum \limits^{N_0} \limits_{i=1} \lambda_i^0 \mathcal L^{2n+1}_{A_i}~~and~~\mu_1=\sum \limits^{N_1} \limits_{i=1} \lambda_i^1 \mathcal L^{2n+1}_{B_i},$$
there exists a unique optimal geodesic plan $\eta\in OptGeo(\mu_0,\mu_1)$ connecting $\mu_0$ and $\mu_1$, and $\Upsilon: =(T_0,T_1)_\sharp\eta\in Opt(\mu_0,\mu_1)$ is the unique optimal transport plan between $\mu_0$ and $\mu_1$ and it is induced by a map $\Pi$, i.e. $\Upsilon=(id\otimes \Pi)_\sharp\mu_0$.

For mapping $\Pi$ we discuss it in two cases:

{\bf Case 1:} If the optimal transport map $\Pi$ is continuous, $\mu_0$ and $\mu_1$ have constant densities, the proof is divided into two parts.

Step 1. From Lemma 3.8 we know $spt(\mu_0)$ is compact, it is possible to find a finite division $\{P^\varepsilon_j\}_{j=1,\ldots,N_3}$ (open set $\{P^\varepsilon_j\}_{j=1,\ldots,N_3}$ mutually disjoint), such that $spt(\mu_0)=\cup_{j=1}^{N_3}P^\varepsilon_j$, and for $j=1,\ldots,N_3$ also has the following properties.

(a) For sets $P^\varepsilon_j,~\Pi(P^\varepsilon_j)$ satisfying: $\mathcal L^{2n+1}(P^\varepsilon_j)>0$, $diam(P^\varepsilon_j)<\varepsilon$ and $diam(\Pi(P^\varepsilon_j))<\varepsilon$;

(b) For $\forall~A_i$, $\exists~j$ such that $P^\varepsilon_j\cap A_i\neq\emptyset$, we take $P^\varepsilon_j\cap A_i\triangleq A_{i,j}^1$ which satisfy the properties $diam(A_{i,j}^1)<\varepsilon$. For the set $A_i\setminus A_{i,j}^1$, if $diam(A_i\setminus A_{i,j}^1)<\varepsilon$, we take $A_i\setminus A_{i,j}^1\triangleq A_{i,j}^2$. If $diam(A_i\setminus A_{i,j}^1)>\varepsilon$, we can divide it into $k$ small pieces according to the property (a). In the same way, for all $A_i,~i=1,\ldots,N_0$, we redivide $A$ into $N_4$ small sets. For convenience we will denote the newly obtained set as $\{P^\varepsilon_j\}$. Similarly, for $\forall~B_l$, $\exists~j$ such that $\Pi(P^\varepsilon_j)\cap B_l\neq\emptyset$, we take $\Pi(P^\varepsilon_j)\cap B_l\triangleq B_{l,j}^1$ and satisfy the properties $diam(B_{l,j}^1)<\varepsilon$. For the set $B_l\setminus B_{l,j}^1$, it is also discussed by situation. In the same way, for all $B_l,~l=1,\ldots,N_1$, we redivide $B$ into $N_5$ small sets and take the newly obtained set to represent the set $\{\Pi(P^\varepsilon_j)\}$. Assume that $\max\{N_3,N_4,N_5\}=L_\varepsilon$, ensuring the properties (a), then we will make appropriate improvements to the segmentation. So, there exists $i(j)$ and $l(j)$ such that $P^\varepsilon_j\subset A_{i(j)}$ and $\Pi(P^\varepsilon_j)\subset B_{l(j)}$, respectively.

By optimal mass transportation, we know that the good transport behavior implies that the unique optimal map $\Pi$ and inverse map $\Phi$:
$$\mu_0(P^\varepsilon_j)\triangleq\mu_0(\Phi\circ\Pi(P^\varepsilon_j))\triangleq\Pi_{\sharp}\mu_0(\Pi(P^\varepsilon_j))\triangleq\mu_1(\Pi(P^\varepsilon_j)).\eqno(4.3)$$
Further, combined with equation (4.3), we can define the measures
$$\mu_0^{\varepsilon,j}:=\mu_0(P^\varepsilon_j)\mathcal L^{2n+1}_{P^\varepsilon_j}~~and~~\mu_1^{\varepsilon,j}:=\mu_1(\Pi(P^\varepsilon_j))\mathcal L^{2n+1}_{\Pi(P^\varepsilon_j)}=\mu_0(P^\varepsilon_j)\mathcal L^{2n+1}_{\Pi(P^\varepsilon_j)}.$$
From the property (b) of this partition, $\mu_0^{\varepsilon,j}$ and $\mu_0|_{P^\varepsilon_j}$ are both measures of constant density with respect to $\mathcal L^{2n+1}$ and with equal mass (the same as $\mu_1|_{\Pi(P^\varepsilon_j)}$ and $\mu_1^{\varepsilon,j}$), we can get
$$\mu_0^{\varepsilon,j}=\mu_0|_{P^\varepsilon_j}~~and~~\mu_0=\sum \limits^{L_\varepsilon} \limits_{j=1}\mu_0|_{P^\varepsilon_j}=\sum \limits^{L_\varepsilon} \limits_{j=1}\mu_0^{\varepsilon,j};$$
$$\Pi_{\sharp}\mu_0^{\varepsilon,j}=\mu_1|_{\Pi(P^\varepsilon_j)}=\mu_1^{\varepsilon,j}~~and~~\mu_1=\sum \limits^{L_\varepsilon} \limits_{j=1}\mu_1|_{\Pi(P^\varepsilon_j)}=\sum \limits^{L_\varepsilon} \limits_{j=1}\mu_1^{\varepsilon,j}.$$
Defining $\eta^\varepsilon_j:=\eta|_{spt(\eta)}\in\mathcal{M}_+(Geo(\mathbb{H}^n))$ and $\eta=\sum \limits^{L_\varepsilon} \limits_{j=1}\eta^\varepsilon_j$. Note that $\bar{\eta}^\varepsilon_j:=\frac{\eta^\varepsilon_j}{\mu_0(P^\varepsilon_j)}$, by (4.3) we have
$$\{\bar{\eta}^\varepsilon_j\}=OptGeo\Big(\frac{(T_0)_\sharp\eta^\varepsilon_j}{\mu_0(P^\varepsilon_j)},\frac{(T_1)_\sharp\eta^\varepsilon_j}{\mu_0(P^\varepsilon_j)}\Big)=OptGeo(\mathcal L^{2n+1}_{P^\varepsilon_j},\mathcal L^{2n+1}_{\Pi(P^\varepsilon_j)}).\eqno(4.4)$$
Thus, $\bar{\mu}_s^{\varepsilon,j}:=(T_s)_\sharp\bar{\eta}^\varepsilon_j$ is the unique Wasserstein geodesic
connecting $\mathcal L^{2n+1}_{P^\varepsilon_j}$ and $\mathcal L^{2n+1}_{\Pi(P^\varepsilon_j)}$, for all $s\in[0,1]$ we note that:
$$spt(\bar{\mu}_s^{\varepsilon,j}): The~support~set~of~measure~\bar{\mu}_s^{\varepsilon,j}.\eqno(4.5)$$
For all $j$ and $s$, from the strong interpolation property, the measure $\bar{\mu}_s^{\varepsilon,j}$ is absolutely continuous and with density $\bar{\rho}_s^{\varepsilon,j}$, and by definition
$$\bar{\rho}_s^{\varepsilon,j}>0~~~~~~\bar{\mu}_s^{\varepsilon,j}-almost~everywhere~on~ spt(\bar{\mu}_s^{\varepsilon,j}).\eqno(4.6)$$
In addition, we can apply the strong Brunn-Minkowski inequality and deduce:
$$\mathcal L^{2n+1}(spt(\bar{\mu}_s^{\varepsilon,j}))^\frac{1}{2n+1}\geq\tau^n_{1-s}(\Theta_j)\mathcal L^{2n+1}(P^\varepsilon_j)^\frac{1}{2n+1}+\tau^n_s(\Theta_j)\mathcal L^{2n+1}(\Pi(P^\varepsilon_j))^\frac{1}{2n+1}, \eqno(4.7)$$
here $\Theta_j:=\Theta_{P^\varepsilon_j,\Pi(P^\varepsilon_j)}$.

Step 2. Next goal is to find a suitable approximately $s$-intermediate point $\tilde{\mu}^\varepsilon_s$ between $\mu_0$ and $\mu_1$. We find this can be achieved by considering a family of measures $\tilde{\mu}^{\varepsilon,j}_s$ supported on the sets $spt(\bar{\mu}_s^{\varepsilon,j})$ and with constant density, then add up them. This would allow us to use (4.7) on each set $P^\varepsilon_j$. For $s\in[0,1]$ and $j=1,\ldots,L_\varepsilon$, we define
$$\tilde{\mu}^{\varepsilon,j}_s:=\mu_0(P^\varepsilon_j)\mathcal L^{2n+1}_{spt(\bar{\mu}_s^{\varepsilon,j})}~~and~~\tilde{\mu}^\varepsilon_s:=\sum \limits^{L_\varepsilon} \limits_{j=1}\tilde{\mu}^{\varepsilon,j}_s.$$
By (4.7), $spt(\bar{\mu}_s^{\varepsilon,j})$ has positive measure. We know $spt(\bar{\mu}_s^{\varepsilon,j})$ is bounded and has finite measure, therefore $\tilde{\mu}^{\varepsilon,j}_s$ is well defined. Since $(\mathbb{H}^n,d,\mathcal L^{2n+1})$ be an essentially non-branching metric measure space, $spt(\bar{\mu}_s^{\varepsilon,i})\cap spt(\bar{\mu}_s^{\varepsilon,j})=\emptyset~(i\neq j)$. Let $\tilde{\rho}^\varepsilon_s$ be the density of $\tilde{\mu}^\varepsilon_s$ with respect to $\mathcal L^{2n+1}$, then
$$\ \ \ \ \ \ Ent_{2n+1}(\tilde{\mu}^\varepsilon_s)=-\int_{\mathbb{H}^n}(\tilde{\rho}^\varepsilon_s)^{1-\frac{1}{2n+1}}d\mathcal L^{2n+1}=-\sum \limits^{L_\varepsilon} \limits_{j=1}\int_{spt(\bar{\mu}_s^{\varepsilon,j})}(\tilde{\rho}^\varepsilon_s)^{1-\frac{1}{2n+1}}d\mathcal L^{2n+1}$$
$$=-\sum \limits^{L_\varepsilon} \limits_{j=1}\mu_0(P^\varepsilon_j)^{1-\frac{1}{2n+1}}\mathcal L^{2n+1}(spt(\bar{\mu}_s^{\varepsilon,j}))^{\frac{1}{2n+1}}.\eqno(4.8)\ \ \ \ \ \ \ \ \ $$
Combination (4.7) and (4.8), $\Upsilon_j:=\Upsilon|_{P^\varepsilon_j\times\Pi(P^\varepsilon_j)}$, gives the following estimate:
$$Ent_{2n+1}(\tilde{\mu}^\varepsilon_s)\leq-\sum \limits^{L_\varepsilon} \limits_{j=1}\mu_0(P^\varepsilon_j)^{1-\frac{1}{2n+1}}\big[\tau^n_{1-s}(\Theta_j)\mathcal L^{2n+1}(P^\varepsilon_j)^\frac{1}{2n+1}+\tau^n_s(\Theta_j)\mathcal L^{2n+1}(\Pi(P^\varepsilon_j))^\frac{1}{2n+1}\big]$$
$$\ =-\sum \limits^{L_\varepsilon} \limits_{j=1}\int[\tau^n_{1-s}(\Theta_j)\rho_0(x)^{-\frac{1}{2n+1}}+\tau^n_{s}(\Theta_j)\rho_1(y)^{-\frac{1}{2n+1}}   ]d\Upsilon_j(x,y)$$
$$\ \ \ \ \ \ \ \ \ \ \leq-\sum \limits^{L_\varepsilon} \limits_{j=1}\int[\tau^n_{1-s}(\theta_x-\varepsilon)\rho_0(x)^{-\frac{1}{2n+1}}+\tau^n_{s}(\theta_x-\varepsilon)\rho_1(y)^{-\frac{1}{2n+1}}   ]d\Upsilon_j(x,y)$$
$$\ \ \ \ \ \ \ \ \ \ \ \ \ \ \ \ \ \ \ =-\int[\tau^n_{1-s}(\theta_x-\varepsilon)\rho_0(x)^{-\frac{1}{2n+1}}+\tau^n_{s}(\theta_x-\varepsilon)\rho_1(y)^{-\frac{1}{2n+1}}   ]d\Upsilon(x,y).\eqno(4.9)$$
Here the first equality follows by the partition in Step 1, $\Upsilon_j$ is concentrated on $P^\varepsilon_j\times\Pi(P^\varepsilon_j)$ and the second inequality is follows from the monotonicity properties of the Heisenberg distortion coefficients.

For general step measures, review the previous notation (4.4), (4.5) and (4.6), for every fixed $t\in[0,1]$ and $j=1,\ldots,L_\varepsilon$, we define the measure
$$\tilde{\eta}^\varepsilon_j:=\frac{\bar{\eta}^\varepsilon_j|_{spt(\bar{\mu}_t^{\varepsilon,j})}}{\mathcal L^{2n+1}(spt(\bar{\mu}_t^{\varepsilon,j}))\bar{\rho}^{\varepsilon,j}_t(T_t(\gamma))}\in OptGeo((T_0)_\sharp\tilde{\eta}_j,(T_1)_\sharp\tilde{\eta}_j).$$
Here $\tilde{\eta}^\varepsilon_j$ is a probability measure. Besides, $\tilde{\eta}^\varepsilon_j$ is concentrated on $Geo(\mathbb{H}^n)$ and $(T_t)_\sharp\tilde{\eta}^\varepsilon_j=\mathcal L^{2n+1}_{spt(\bar{\mu}_t^{\varepsilon,j})}$.
In addition, the measures defined by $\nu^{\varepsilon,j}_0:=(T_0)_\sharp\tilde{\eta}_j^\varepsilon$ and $\nu^{\varepsilon,j}_1:=(T_1)_\sharp\tilde{\eta}_j^\varepsilon$ are concentrated on $P^\varepsilon_j$ and $\Pi(P^\varepsilon_j)$ respectively. From the Wasserstein distance, for $\omega_1,\omega_2\in\mathcal{P}_2(\mathbb{H}^n)$ it holds
$$W_2(\omega_1,\omega_2)\leq diam(spt(\omega_1)\cup spt(\omega_2)),\eqno(4.10)$$
Applying property (a) and triangular inequality, for $j=1,\ldots,L_\varepsilon$, we can get
$$W_2(\mathcal L^{2n+1}_{P^\varepsilon_j},\mathcal L^{2n+1}_{spt(\bar{\mu}_t^{\varepsilon,j})})\leq W_2(\mathcal L^{2n+1}_{P^\varepsilon_j},\nu^{\varepsilon,j}_0)+W_2(\nu^{\varepsilon,j}_0,\mathcal L^{2n+1}_{spt(\bar{\mu}_t^{\varepsilon,j})})$$
$$\ \ \ \ \leq\varepsilon+tW_2(\nu^{\varepsilon,j}_0,\nu^{\varepsilon,j}_1)$$
$$ \ \ \ \ \ \ \ \ \  \ \ \ \ \ \leq3\varepsilon+tW_2(\mathcal L^{2n+1}_{P^\varepsilon_j},\mathcal L^{2n+1}_{\Pi(P^\varepsilon_j)}),\eqno(4.11)$$
Similarly, inequality (4.10) is used to obtain
$$W_2(\mathcal L^{2n+1}_{spt(\bar{\mu}_t^{\varepsilon,j})},\mathcal L^{2n+1}_{\Pi(P^\varepsilon_j))})\leq3\varepsilon+(1-t)W_2(\mathcal L^{2n+1}_{P^\varepsilon_j},\mathcal L^{2n+1}_{\Pi(P^\varepsilon_j)}).$$
Now recalling that $\mu_0=\sum \limits^{L_\varepsilon} \limits_{j=1}\mu_0(P^\varepsilon_j)\mathcal L^{2n+1}_{P^\varepsilon_j}$ and $\tilde{\mu}^\varepsilon_t=\sum \limits^{L_\varepsilon} \limits_{j=1}\mu_0(P^\varepsilon_j)\mathcal L^{2n+1}_{spt(\bar{\mu}_t^{\varepsilon,j})}$, the convexity of $W_2^2$ gives us
$$W_2^2(\mu_0,\tilde{\mu}^\varepsilon_t)\leq\sum \limits^{L_\varepsilon} \limits_{j=1}\mu_0(P^\varepsilon_j)W_2^2(\mathcal L^{2n+1}_{P^\varepsilon_j},\mathcal L^{2n+1}_{spt(\bar{\mu}_t^{\varepsilon,j})}).$$
Because $\Pi$ mapping is optimal, we have
$$\sum \limits^{L_\varepsilon} \limits_{j=1}\mu_0(P^\varepsilon_j)W_2^2(\mathcal L^{2n+1}_{P^\varepsilon_j},\mathcal L^{2n+1}_{\Pi(P^\varepsilon_j)})=W_2^2(\mu_0,\mu_1).$$
For $t\in[0,1]$, summing (4.11) on  all $j$, we obtain
$$W_2^2(\mu_0,\tilde{\mu}^\varepsilon_t)\leq\sum \limits^{L_\varepsilon} \limits_{j=1}\mu_0(P^\varepsilon_j)\Big(3\varepsilon+tW_2(\mathcal L^{2n+1}_{P^\varepsilon_j},\mathcal L^{2n+1}_{\Pi(P^\varepsilon_j)})\Big)^2\ \ \ \ \ \ \ \ \ \ \ \ \ \ \ \ \ \ \ \ \ \ \ \ \ \ \ \ \ \ \ \ \ \ \ \ \ \ \ \ \ \ \ \ \ $$
$$ \ \ \ \ \ \ \ \ \ \ \ \ \ \ \ =9\varepsilon^2+6\varepsilon t\sum \limits^{L_\varepsilon} \limits_{j=1}\mu_0(P^\varepsilon_j)W_2(\mathcal L^{2n+1}_{P^\varepsilon_j},\mathcal L^{2n+1}_{\Pi(P^\varepsilon_j)})+t^2\sum \limits^{L_\varepsilon} \limits_{j=1}\mu_0(P^\varepsilon_j)W_2^2(\mathcal L^{2n+1}_{P^\varepsilon_j},\mathcal L^{2n+1}_{\Pi(P^\varepsilon_j)})$$
$$\leq9\varepsilon^2+6\varepsilon t\cdot diam(spt(\mu_0)\cup spt(\mu_1))+t^2W_2^2(\mu_0,\mu_1). \ \ \ \ \ \ \ \ \ \ \ \ \ \ \eqno(4.12)$$
Analogously, we also have
$$W_2^2(\tilde{\mu}^\varepsilon_t,\mu_1)\leq9\varepsilon^2+6\varepsilon (1-t)\cdot diam(spt(\mu_0)\cup spt(\mu_1))+(1-t)^2W_2^2(\mu_0,\mu_1).\eqno(4.13)$$
From Lemma 3.8, for $s\in[0,1]$, $\mu_0$ and $\mu_1$ have bounded support and the space $(\mathbb{H}^n,d,\mathcal L^{2n+1})$ is proper, all measures family $\{\tilde{\mu}^\varepsilon_s\}_{\varepsilon>0}$ are concentrated on a common compact set. So we can find a sequence $\{\varepsilon_m\}_{m\in\mathbb{N}}$ converging to $0$ such that
$$\tilde{\mu}^{\varepsilon_m}_s \stackrel{^{W_2}}\rightarrow \mu_s\in\mathcal{P}_2(\mathbb{H}^n)~~~~as~~~m\rightarrow\infty.$$
As $m\rightarrow\infty$, from (4.12) and (4.13) we get
$$W_2(\mu_0,\mu_s)\leq sW_2(\mu_0,\mu_1)~~and~~W_2(\mu_s, \mu_1)\leq (1-s)W_2(\mu_0,\mu_1).$$
Combining this two inequalities with the triangular inequality we deduce that
$$W_2(\mu_0,\mu_s)= sW_2(\mu_0,\mu_1)~~and~~W_2(\mu_s, \mu_1)= (1-s)W_2(\mu_0,\mu_1),$$
which means that $\mu_s$ is the unique $s$-midpoint of $\mu_0$ and $\mu_1$, and
$$\tilde{\mu}^{\varepsilon}_s \stackrel{^{W_2}}\rightarrow \mu_s\in\mathcal{P}_2(\mathbb{H}^n)~~~~as~~~\varepsilon\rightarrow 0.$$
For all $s\in[0,1]$, we deduce that the curve $t\mapsto\mu_s$ is the unique Wasserstein geodesic connecting $\mu_0$ and $\mu_1$. And then take the limit of inequality (4.9),
$$Ent_{2n+1}(\mu_s)\leq-\int\Big[\tau^n_{1-s}(\theta_x)\rho_0(x)^{1-\frac{1}{2n+1}}+\tau^n_s(\theta_x)\rho_1(y)^{1-\frac{1}{2n+1}}
\Big]d\Upsilon(x,y).$$
This leads to the desired inequality. So we have proved that if the optimal mapping between two step measures with bounded support is continuous, the unique geodesic connecting them satisfie (3.1).

{\bf Case 2:} If $\Pi$ is not continuous, from Lemma 3.8 we know $\mathcal L^{2n+1}$ is Radon, apply Lusin theorem: for $\forall~\varepsilon>0$, there exists a compact set $A_\varepsilon\subset spt(\mu_0)$, such that the mapping $\Pi$ is continuous on $A_\varepsilon$ and $\mu_0(spt(\mu_0)\setminus{A_\varepsilon})<\varepsilon$.
Then, define the unique optimal geodesic plan
$$\eta^\varepsilon:=\frac{\eta|_{spt(\eta)}}{\mu_0(A_\varepsilon)}\in\mathcal{P}(Geo(\mathbb{H}^n))~~and~~\mu_0^\varepsilon:=(T_0)_\sharp\eta,~\mu_1^\varepsilon:=(T_1)_\sharp\eta\in\mathcal{P}^{ac}(\mathbb{H}^n,\mathcal L^{2n+1}).$$
Next, exploiting the good transport behaviour as done in the first part, we deduce the set $B_\varepsilon:\Pi(A_\varepsilon)\subset spt(\mu_1)$ such that $\mu_1^\varepsilon=\mu_1|_{B_\varepsilon}/\mu_1(B_\varepsilon)$, while $\mu_0^\varepsilon=\mu_0|_{A_\varepsilon}/\mu_0(A_\varepsilon)$. In particular, $\mu_0^\varepsilon$ and $\mu_1^\varepsilon$ are step measures with bounded support, $\Pi|_{A_\varepsilon}$ is the optimal map continuous. Therefore inequality (3.1) holds for $\mu_0^\varepsilon$, $\mu_1^\varepsilon$ and $\eta^\varepsilon$. For $A_\varepsilon$ and $B_\varepsilon$, we use the method of proposition 3.6, as $\varepsilon\rightarrow0$, (3.1) holds for $\mu_0$, $\mu_1$ and $\eta$. This concludes the proof. 

(ii) For each pair of bounded measure $\mu,~\nu\in\mathcal{P}^{ac}(\mathbb{H}^n,\mathcal L^{2n+1})$, there exist two approximating sequences $\{\mu_m\}_{m\in \mathbb{N}}$, $\{\nu_m\}_{m\in \mathbb{N}}$  satisfying the requirements of Lemma 3.4. For every $m\in \mathbb{N}$, $\eta_m$ is the unique optimal geodesic plan in $OptGeo(\mu_m,\nu_m)$, let $\Upsilon_m:=(T_0,T_1)_\sharp\eta_m$ and $\mu_s^m:=(T_s)_\sharp\eta_m$. For every $m\in\mathbb{N}$, $\mu_m$ and $\nu_m$ are step measures with bounded support, then according to (i), for all $s\in[0,1]$ there is
$$Ent_{2n+1}(\mu_s^m)\leq\mathcal{F}^n_s(\Upsilon_m|\mathcal L^{2n+1}). \eqno(4.14)$$
Significantly,  $\{\mu_m\}_{m\in \mathbb{N}}\bigcup\{\mu\}$ and $\{\nu_m\}_{m\in \mathbb{N}}\bigcup\{\nu\}$ are uniformly bounded and the space $(\mathbb{H}^n,d, \mathcal L^{2n+1})$ is proper for every $s\in[0,1]$, all the measures in the family $\{\mu_s^m\}_{m\in \mathbb{N}}$ are concentrated on the same compact set. In particular, for every $s\in[0,1]$, the family $\{\mu_s^m\}_{m\in \mathbb{N}}$ is $W_2$-precompact, then
$$\mu_s^m \stackrel{^{W_2}}\rightarrow \mu_s\in\mathcal{P}_2(\mathbb{H}^n)~~as~~m\rightarrow\infty.$$
$\mu_s^m$ is a $s$-midpoint between $\mu_m$ and $\nu_m$, $\mu_m\stackrel{^{W_2}}\rightarrow \mu$ and $\nu_m\stackrel{^{W_2}}\rightarrow \nu$, thus $\mu_s$ is the unique $s$-midpoint between $\mu$ and $\nu$.
Then, due to the lower semicontinuity of entropy functional and lemma 3.5, we can take the limit of $m\rightarrow\infty$ for formula (4.14), and obtain
$$Ent_{2n+1}(\mu_s)\leq\mathcal{F}^n_s(\Upsilon|\mathcal L^{2n+1})~~~for~all~s\in[0,1]$$
which is inequality (3.1) and conclude the proof.
\hfill ${\square}$

\noindent{\bf Remark 4.4}.~~If $A$ and $B$ are translated sets, from the define of Heisenberg distortion coefficients and $\Theta_{A,B}\neq 0$, we know the equality of inequality (3.1) do not hold.

\end{document}